\documentclass{article}
\usepackage{graphicx}
\usepackage[T1]{fontenc}
\usepackage[english]{babel}
\usepackage{graphicx}
\usepackage{amssymb, amsmath, amsfonts}
\usepackage{dsfont}
\usepackage{vmargin}
\usepackage{makeidx}

\everymath{\displaystyle}

\usepackage{hyperref}
\hypersetup{
	colorlinks=true,
	breaklinks=true,
	urlcolor=blue,
	linkcolor=red,
	bookmarksopen=true,
	pdftitle = Stochastic Cahn-Hilliard equation with double singular nonlinearities and two reflections,
	pdfauthor = Ludovic Gouden\`ege,
	pdfsubject = Stochastic Cahn-Hilliard equation with double singular nonlinearities and two reflections,
	pdfkeywords = {Cahn-Hilliard, stochastic partial differential equations, integration by parts formulae, reflection measures, invariant measures, singular nonlinearity, double, two}
}

\numberwithin{equation}{section} 

\newcommand{\R}{\mathbb{R}} 
\newcommand{\B}{\mathbf{B}} 
\newcommand{\N}{\mathbb{N}} 
\newcommand{\E}{\mathbb{E}} 
\newcommand{\F}{\mathcal{F}} 
\newcommand{\1}{\mathds{1}} 
\newcommand{\dd}{\mathrm{d}} 
\newcommand{\dds}{\mathrm{ds}}
\newcommand{\ddt}{\mathrm{dt}}
\newcommand{\ddtheta}{\mathrm{d\theta}}
\newcommand{\ddx}{\mathrm{dx}}
\newcommand{\ddy}{\mathrm{dy}}

\newcommand{\BD}{\begin{displaymath}}
\newcommand{\ED}{\end{displaymath}}
\newcommand{\BEA}{\begin{eqnarray}}
\newcommand{\EEA}{\end{eqnarray}}
\newcommand{\BEAS}{\begin{eqnarray*}}
\newcommand{\EEAS}{\end{eqnarray*}}
\newcommand{\BE}{\begin{equation}}
\newcommand{\EE}{\end{equation}}
\newcommand{\BES}{\begin{equation*}}
\newcommand{\EES}{\end{equation*}}

\newcommand{	\BCB}{\color{black}}

\newcommand{\Proof}{\noindent {\bf{\underline{Proof} : }}}
\newcommand{\EProof}{\begin{flushright}$\Box$\end{flushright}}

\newcommand{\Step[1]}{\indent \textbf{Step {#1}.}}
\newcommand{\EStep[1]}{\begin{flushright}$\Box$\end{flushright}}

\def\ds{\displaystyle}
\def\P{\mathbb P}
\def\Tr{\mbox{Tr}}

\newtheorem{Def}{Definition}[section]
\newtheorem{Coro}{Corollary}[section]
\newtheorem{Le}{Lemma}[section]
\newtheorem{Th}{Theorem}[section]
\newtheorem{Prop}{Proposition}[section]
\newtheorem{Rem}{Remark}[section]

\title{Stochastic Cahn-Hilliard equation with double singular nonlinearities and two reflections}
\author{Arnaud Debussche, Ludovic Gouden\`ege\\
{\small ENS Cachan - Antenne de Bretagne, Campus de Ker-Lann, 35170 Bruz, France}\\
{\small arnaud.debussche@bretagne.ens-cachan.fr, ludovic.goudenege@bretagne.ens-cachan.fr}}
\date{}

\begin{document}

\maketitle
\footnotetext{\!\!\!\!\!\!\!\!\!\!\!AMS 2000 subject classifications. 60H15, 60H07, 37L40.}

\section*{R\'esum\'e}
On consid\`ere une \'equation aux d\'eriv\'ees partielles stochastique poss\'edant deux non-lin\'earit\'es de type logarithmique, avec deux r\'eflexions en $1$ et $-1$ sous la contrainte de conservation de masse. L'\'equation, dirig\'ee par un bruit blanc en espace et en temps, contient un double Laplacien. L'absence de principe de maximum pour le double Laplacien pose des difficult\'es pour l'utilisation d'une m\'ethode classique de p\'enalisation, pour laquelle une importante propri\'et\'e de monotonie est utilis\'ee. Etant inspir\'e par les travaux de Debussche, Gouden\`ege et Zambotti, on d\'emontre l'existence et l'unicit\'e de solutions pour des donn\'ees initiales entre $-1$ et $1$. Enfin, on d\'emontre que l'unique mesure invariante est ergodique, et on \'enonce un r\'esultat de m\'elange exponentiel.
\section*{Abstract}
We consider a stochastic partial differential equation with two logarithmic nonlinearities, with two reflections at $1$ and $-1$ and with a constraint of conservation of the space average. The equation, driven by the derivative in space of a space-time white noise, contains a bi-Laplacian in the drift. The lack of the maximum principle for the bi-Laplacian generates difficulties for the classical penalization method, which uses a crucial monotonicity property. Being inspired by the works of Debussche, Gouden\`ege and Zambotti, we obtain existence and uniqueness of solution for initial conditions in the interval $(-1,1)$. Finally, we prove that the unique invariant measure is ergodic, and we give a result of exponential mixing.
\section*{Introduction and main results}
The Cahn-Hilliard-Cook equation is a model to describe phase separation in a binary alloy (see \cite{C}, \cite{CH1} and \cite{CH2}) in the presence of thermal fluctuations (see \cite{COOK} and \cite{LANGER}). It takes the form:
\BE
\label{e0.1}
\left\{\begin{array}{ll}
\partial_{t} u = - \frac{1}{2}\Delta\left(\Delta u - \psi(u)\right)+\dot{\xi},&\text{ on } \Omega\subset \R^n,\\
\\
\nabla u \cdot \nu = 0 = \nabla (\Delta u-\psi(u)) \cdot \nu, &\text{ on } \partial\Omega,\\ 
\end{array}\right.
\EE
where $t$ denotes the time variable and $\Delta$ is the Laplace operator. Also $u\in [-1,1]$ represents the ratio between the two species and the noise term $\dot{\xi}$ accounts for the thermal fluctuations. The nonlinear term $\psi$ has the double-logarithmic form:
\BE
\label{e0.1bis}
\psi : u \mapsto \frac{\theta}{2} \ln\left(\frac{1+u}{1-u}\right) - \theta_{c} u,
\EE
where $\theta$ and $\theta_{c}$ are temperatures with $\theta < \theta_{c}$.

The study of this equation presents several difficulties. First, the singularities at $\pm 1$ have to be treated
carefully. Also, since it is a fourth order equation, no comparison principle holds.

The deterministic equation where $\psi$ is replaced by a polynomial function have first been studied (see \cite{CH1}, \cite{LANGER} and \cite{MR763473}). 
Then non smooth $\psi$ have been considered (see \cite{MR1123143} and \cite{MR1327930}).

Phase separation have been analysed thanks to this model: see for example the survey \cite{MR1657208}, and the references therein, or others recent results on spinodal decomposition and nucleation in \cite{MR1232163, MR2342011, MR1214868, MR1637817, MR1753703, MR1712442, MR1763320, MR2048517}. 

In the case of a polynomial nonlinearity, some results  have been obtained in the stochastic case (see \cite{BLMAWA1, BLMAWA2, MR1867082, MR1897915, MR1359472, MR1111627}). 

Note that the solutions of 
the equation with polynomial nonlinearity do not remain in $[-1,1]$ in general, and their physical interpretation
is not clear.

To our knowledge, the case of the logarithmic nonlinearity in the presence of noise have never been studied. The presence of noise has a strong effect and equation \eqref{e0.1} cannot have a solution. Indeed,
a solution should remain in $[-1,1]$ which is impossible with an additive noise. Two reflection measures 
have to be added to the model to remedy this problem. 
In this article, we propose to study:
\BE\label{Eq:0.2}
\left\{\begin{array}{ll}
\displaystyle{\partial_{t} u = - \frac{1}{2}\Delta\Big(\Delta u -\psi(u) + \eta_{-} - \eta_{+} \Big) + \dot{\xi}},&\text{ with } \theta \in [0,1]=\Omega,\\
&\\ 
\nabla u \cdot \nu = 0 = \nabla (\Delta u) \cdot \nu, &\text{ on } \partial\Omega,
\end{array}\right.
\EE
where the measures are subject to the contact conditions almost surely:
\BE\label{Eq:0.3}
\int (1+u) \mathrm{d}\eta_{-} = \int (1-u) \mathrm{d}\eta_{+} = 0.                                                                                                                                                                                                                                                                                              
\EE

The stochastic heat equation with reflection,  {\it i.e.} when the fourth order operator is replaced 
by the Laplace operator, is a model for the evolution of random interfaces near a hard wall. 
It has been extensively studied in the literature (see \cite{MR2257651}, \cite{MR1835843}, \cite{MR1872740},  \cite{MR1172940} \cite{MR1921014}, \cite{MR1959795} and \cite{MR2128236}). 
Essential tools in these articles are  the comparison principle and the fact that the underlying Dirichlet form 
is symmetric so that the invariant measure is known explicitely.

In our case, we consider a noise which is obtained as the space derivative of the space-time
white noise. In other words, the noise is the time derivative of a cylindrical Wiener process in 
$H^{-1}(0,1)$. This is physically reasonable since the Cahn-Hilliard equation can be interpreted 
as a gradient system in this space. With such noise, the system is still symmetric and the invariant measure
is known explicitely. As in the second order case, we use this fact in an essential way.

However, as already mentioned, no comparison principle holds and new techniques have to be developed.
The equation \eqref{Eq:0.2} has been studied with a single reflection and when no nonlinear term is taken into account in \cite{MR2349572}. The reflection is introduced to enforce positivity of the 
solution. Various techniques have been introduced to  overcome this lack of comparison principle. 
Moreover, as in the second order case, an integration by part formula for the invariant measure has been 
derived. Then, in \cite{GOUDENEGE}, a singular nonlinearity of the form $u^{-\alpha}$ or $\ln u$ 
have been considered. Existence and uniqueness of solutions have been obtained and using the 
integration by parts formula as in \cite{MR1959795}, it has been proved that the reflection measure 
vanishes if and only if $\alpha \ge 3$. In particular, for a logarithmic nonlinearity, the reflection is active. 

Here, we consider the original Cahn-Hilliard-Cook model \eqref{e0.1} with the double-logarithmic nonlinear 
term  \eqref{e0.1bis}. The noise is as in the above mentioned articles and we still have an explicit 
invariant measure.  Our method mixes ideas from \cite{MR2349572}, \cite{GOUDENEGE} and \cite{MR1959795}. Additional difficulties are overcome, the main one being to understand how to deal with the nonlinear term. Indeed, in \cite{GOUDENEGE}, the positivity of the nonlinear term was essential.
We overcome this difficulty thanks to a delicate a priori estimate. 
Our main results state that equations \eqref{Eq:0.2}, \eqref{Eq:0.3} together with an initial condition have a unique solution (see Proposition \ref{Prop:2.3} and Theorem \ref{Th:Existence}). As in \cite{MR2349572}, it 
is constructed thanks to the gradient structure of \eqref{Eq:0.2} and strong Feller property. Moreover, we prove that this solution is the limit of the solution of the Cahn-Hilliard-Cook equation with polynomial 
nonlinearity without reflections. This justifies the use of the polynomial models.
We also prove that the invariant measure is unique and ergodic. Such property is very easy to obtain if 
$\theta_c$ is small (see  \cite{MR2349572}) or in the polynomial case (see \cite{MR1359472}). Finally, a stronger result of exponential mixing is given in the last Theorem \ref{Th:Mixing}. It is based on coupling and arguments developped by Odasso in \cite{MR2276442}.

In future studies, we shall generalize the integration by part formula obtain in \cite{MR2349572} to 
prove that the reflection measure does not vanish. The presence of two reflection measures introduces
additional difficulties. In the second order case, this has been studied in \cite{MR2278459}.

\section{Preliminaries}\label{S:1}

We denote by $\langle\cdot,\cdot\rangle$ the scalar product in $L^2(0,1)$;  $A$ is the realization in $L^2(0,1)$ of the Laplace operator with Neumann boundary condition, i.e.: 
\BD
D(A) = \text{ Domain of }A  = \{ h \in W^{2,2}(0,1) : h'(0) = h'(1) = 0\}
\ED
where we use $W^{n,p}$ and $||.||_{W^{n,p}}$ to denote the Sobolev space $W^{n,p}(0,1)$ and its associated norm. 
Remark that $A$ is self-adjoint on $L^2(0,1)$ and we have a complete orthonormal system of eigenvectors $(e_{i})_{i\in\N}$ in $L^2(0,1)$ for the eigenvalues $(\lambda_{i})_{i\in\N}$. We denote by $\bar{h}$ the mean of $h \in L^2(0,1)$: 
\BD
\bar{h} = \int_{0}^1h(\theta)\ddtheta.
\ED
We remark that $A$ is invertible on the space of functions with $0$ average. In general, we define
$(-A)^{-1}h=(-A)^{-1}(h-\bar h)+\bar h$.

For $\gamma \in\R$, we define $(-A)^{\gamma}$ by classical interpolation. We set 
$V_{\gamma}:=D((-A)^{\gamma/2})$. 
It is endowed with the classical seminorm and norm :
\BD
|h|_{\gamma} =  \left(\sum_{i=1}^{+\infty} (-\lambda_{i})^{\gamma} h_{i}^2\right)^{1/2},\; 
\|h\|_{\gamma} =  \left(|h|_{\gamma}^2 + \bar{h}^2\right)^{1/2},\mbox{ for } h=\sum_{i\in\N} h_i e_i.
\ED
$|\cdot|_\gamma$ is associated to the scalar product $( \cdot,\cdot)_{\gamma}$.
To lighten notations, we set $(\cdot,\cdot) := (\cdot,\cdot)_{-1}$  and $H:=V_{-1}$. The average can be defined in any $V_\gamma$ by $\bar h= (h,e_0)$. It plays an important role and we often work with functions with a fixed average $c\in \R$. We define $H_c=\{h\in H : \bar{h}=c\}$ for all $c\in \R$.

We use the following regularization operators:
$$
Q_Nx= \frac1N\sum_{n=0}^N\sum_{i=0}^n(x,e_i)e_i.
$$
It is defined on $L^2(0,1)$ and can extended to any $V_\gamma$. Clearly
$Q_Nx$ converges to $x$ in $V_\gamma$ if $x\in V_\gamma$. Moreover, it is well known that
if $x\in C([0,1];\R)$, then the converges holds in $C([0,1];\R)$. Note also that $Q_N$ is self-adjoint
in $V_\gamma$ and commutes with $A$.

The covariance operator of the noise is  the operator $B$ defined by
\BD
B = \frac{\partial}{\partial \theta}, D(B) = W_0^{1,2}(0,1).
\ED
Note that 
\BD
B^*=-\frac{\partial}{\partial \theta},\, D(B^*) = W^{1,2}(0,1), BB^*=-A.
\ED
We denote by $\mathcal{B}_b(H_{c})$ the space of all Borel bounded functions 
on $H_c$. We set $O_{s,t} := [s,t] \times [0,1]$ for $s,t \in [0,T]$ with $s<t$ and $T>0$, and $O_{t} = O_{0,t}$ for $0\leq t\leq T$. Given a measure $\zeta$ on $O_{s,t}$ and a continuous function $v$ on $O_{s,t}$, we write
\BD
\big\langle v,\zeta\big\rangle_{O_{s,t}} := \int_{O_{s,t}} v\ \dd\zeta.
\ED

For $\lambda \in \R$, we define:
\BE\label{Eq:0.4}                                                                                                                                                                                                                                                                                                      
f(x):=\left\{\begin{array}{lr}
+\infty,& \text{ for all } x \leq -1,\\
\\
\ln \left(\frac{1-x}{1+x}\right) +Ê\lambda x,& \text{ for all } x \in (-1,1),\\
\\
-\infty,& \text{ for all } x \geq 1,
\end{array}\right.                                                                                                                                                                                                                                                                                                      
\EE
and the following antiderivative $F$ of $-f$:
\BD
F(x)= (1+x) \ln(1+x)+(1-x) \ln(1-x)-\frac{\lambda}2 x^2,\text{ for all }x \in (-1,1).
\ED
With these notations, we rewrite \eqref{Eq:0.2} in the abstract form:
\BE\label{Eq:1.1}
\left\{\begin{array}{l}
\dd X=-\frac{1}{2}A\left(AX+f(X)+\eta_{-}-\eta_{+}\right)\ddt + B\dd W,\\
\\
\langle (1+X) ,\eta_{-}\rangle_{O_{T}} = \langle (1-X) ,\eta_{+}\rangle_{O_{T}} =0,\\
\\
X(0,x)=x \text{ for } x \in V_{-1},
\end{array}\right.
\EE
where $W$ is a cylindrical Wiener process on $L^2(0,1)$.
\begin{Def}\label{De:1}
Let $x \in \mathcal{C}([0,1];[-1,1])$. We say that $\left(\left(X(t,x)\right)_{t\in[0,T]},\eta_{+},\eta_{-},W\right)$, defined on a filtered complete probability space $\left(\Omega, \P, \F, (\F_t)_{t\in [0,T]}\right)$, is a weak solution to \eqref{Eq:0.2} on $[0,T]$ for the initial condition $x$ if:
\begin{enumerate}
\item[(a)] a.s. $X \in \mathcal{C}\left((0,T]\times[0,1];[-1,1]\right) \cap \mathcal{C}([0,T];H)$ and $X(0,x) = x$,
\item[(b)] a.s. $\eta_{\pm}$ are two positive measures on $(0,T]\times [0,1]$, such that $\eta_{\pm}(O_{\delta,T}) <+\infty$ for all $\delta \in (0,T]$,
\item[(c)] $W$ is a cylindrical Wiener process on $L^2(0,1)$,
\item[(d)] the process $\left(X(\cdot,x),W\right)$ is $(\F_t)$-adapted,
\item[(e)] a.s. $f(X(\cdot,x)) \in L^1(O_{T})$,
\item[(f)] for all $h \in D(A^2)$ and for all $0 < \delta \leq t \leq T$ :
\BD
\begin{array}{rcl}
\langle X(t,x),h\rangle &=& \langle X(\delta,x),h\rangle - \frac{1}{2}\int_{\delta}^t\langle X(s,x),A^2h\rangle \dds-\frac{1}{2}\int_{\delta}^t\langle Ah,f(X(s,x))\rangle \dds\\
\\
&&- \frac{1}{2}\big\langle Ah,\eta_{+}\big\rangle_{O_{\delta,t}} + \frac{1}{2}\big\langle Ah,\eta_{-}\big\rangle_{O_{\delta,t}} - \int_{\delta}^t\langle Bh,\dd W\rangle, \quad a.s.,
\end{array}
\ED
\item[(g)] a.s. the contact properties hold : \\
$supp(\eta_{-}) \subset \{(t,\theta)\in O_{T} / X(t,x)(\theta)=-1\}$ and $supp(\eta_{+}) \subset \{(t,\theta)\in O_{T} / X(t,x)(\theta)=1\}$, that is, 
\BD
\big\langle (1+X),\eta_{-}\big\rangle_{O_{T}}=\big\langle (1-X),\eta_{+}\big\rangle_{O_{T}}= 0.
\ED
\end{enumerate}
Finally, a weak solution $(X,\eta_{+},\eta_{-},W)$ is a strong solution if the process $t \mapsto X(t,x)$ is adapted to the filtration $t \mapsto \sigma(W(s,.),s\in [0,t])$
\end{Def}
\BCB
\begin{Rem} In (f), the only term where we use the function $f$ is well defined. Indeed, by (e) we have $f(X(\cdot,x)) \in L^1(O_{T})$ and by Sobolev embedding $Ah \in D(A) \subset L^\infty(O_{T})$. Hence the notation $\langle \cdot, \cdot \rangle$ should be interpreted as a duality between $L^\infty$ and $L^1$.
\end{Rem}

The solution of the linear equation with initial data $x \in H$ is given by
\BES
Z(t,x) = e^{-tA^2/2}x + \int_0^te^{-(t-s)A^2/2}B\dd W_s.
\EES
As easily seen this process is in $\mathcal{C}([0,+\infty[;H)$ (see \cite{MR1207136}). In particular, the mean of $Z$ is constant and the law of the process $Z(t,x)$ is the Gaussian measure:
\BD
Z(t,x)\sim\mathcal{N}\big(e^{-tA^2/2}x,Q_t\big),\;
Q_t = \int_0^te^{-sA^2/2}BB^*e^{-sA^2/2}\dds = (-A)^{-1}(I - e^{-tA^2}).
\ED
If we let $t\rightarrow +\infty$, the law of $Z(t,x)$ converges to the Gaussian measure on $L^2$:
\BD
\mu_c := \mathcal{N}(ce_0,(-A)^{-1}), \text{ where } c=\bar{x}.
\ED
Notice that $\mu_c$ is concentrated on $H_c\cap C([0,T])$.

In order to solve equation \eqref{Eq:1.1}, we use polynomial approximations of this equation. We denote by $\{f^n\}_{n \in \N}$ the sequence of polynomial functions which converges to the function $f$ on $(-1,1)$, defined for $n \in \N$ by: 
\BD
f^n(x)=-2\sum_{k=0}^{n}\frac{x^{2k+1}}{(2k+1)} + \lambda x, \text{ for all } x \in \R.
\ED
We use the following antiderivative $F^n$ of $-f^n$ defined by:
\BD
F^n(x)=2\sum_{k=0}^{n}\frac{x^{2k+2}}{(2k+2)(2k+1)}-\frac\lambda2 x^2, \text{ for all } x \in \R.
\ED

\noindent Then for $n \in \N$, we study for the following polynomial approximation of \eqref{Eq:1.1} with an initial condition $x \in H$: 
\BE\label{Eq:1.2}
\left\{\begin{array}{l}
\dd X^n+\frac{1}{2}(A^2X^n+Af^n(X^n))\ddt = B\dd W,\\
\\
X^n(0,x)=x.
\end{array}\right.
\EE
This equation has been studied in \cite{MR1359472} in the case $B=I$. The results generalize immediately and  it can be proved that for any $x\in H$, there exists a unique solution $X^n(\cdot,x)$ $a.s.$ in 
$\mathcal{C}([0,T];H)\cap  L^{2n+2}((0,T)\times (0,1))$. 
It is a solution in the mild or weak sense. Moreover the average of $X^n(t,x)$ does not depend on $t$.

For each $c\in \R$, \eqref{Eq:1.1} defines a transition  semigroup $(P^{n,c}_t)_{t\ge 0}$:
\BES
P^{n,c}_t\phi(x) = \E[\phi(X^{n}(t,x)], \; t \geq 0, x \in H_{c}, \; \phi \in \mathcal{B}_b(H_{c}),\; n \in \N^*.
\EES
Existence of an invariant measure can be proved as in \cite{MR1359472}.

Using Galerkin approximation and Bismut-Elworthy-Li formula, it can be seen that $(P^{n,c}_t)_{t\ge 0}$ is Strong Feller. 
More precisely,  for all $\phi \in \mathcal{B}_{b}(H_{c})$,  $n\in\N$ and $t >0$: 
\BE\label{Eq:1.3}
| P^{n,c}_{t}\phi(x) - P^{n,c}_{t}\phi(y)| \leq \frac{2e^{\lambda^2t/4}}{\lambda\sqrt{t}}\|\phi\|_\infty |x-y|_{-1}, \quad\text{for all } x,y \in H_{c}.
\EE
Irreducibility follows from a control argument. By Doob Theorem we deduce that there 
exists an unique and ergodic invariant measure $\nu^n_c$.

It is classical that equation \eqref{Eq:1.2} is a gradient system in $H_c$  and can be rewritten as:
\BE
\left\{\begin{array}{l}
\dd X^n+\frac{1}{2}A(AX^n-\nabla U^n(X^n))\ddt = B\dd W,\\
\\
X^n(0,x)= x \in L^2(0,1),
\end{array}\right.
\EE
where $\nabla$ denotes the gradient in the Hilbert space $L^2(0,1)$, and: 
\BD
U^n(x):= \int_{0}^{1}F^n(x(\theta))\ddtheta, \quad x \in L^2(0,1).
\ED
The measure $\nu^n_c$ is therefore given by:
\BD
\nu_c^n(\ddx)=\frac{1}{Z_c^n}\exp(-U^n(x))\mu_c(\ddx),
\ED
where $Z_c^n$ is a normalization constant. 

We prove in section \ref{S:2} that, for $c\in (-1,1)$, the sequence $(\nu_c^n)_{n\in\N}$ converges to the measure
$$
\nu_c(\ddx)=\frac{1}{Z_c}\exp(-U(x))\1_{x\in K}\mu_c(\ddx),
$$
where 
$$
U(x):= \int_{0}^{1}F(x(\theta))\ddtheta, \quad x \in L^2(0,1).
$$
and 
$$
K=\{x\in L^2 :\; 1\ge x \ge -1\}.
$$
In section \ref{S:2}, we prove the following result.

\begin{Th}\label{Th:Existence} Let $c\in (-1,1)$.
Let $x\in K$ such that $\bar x=c$, then there exists a continuous process denoted $(X(t,x))_{t\geq 0}$ and two nonnegative measures $\eta^x_{+}$ and $\eta^x_{-}$ such that 
$\left(\left(X(t,x)\right)_{t\geq0},\eta^x_{+},\eta^x_{-},W\right)$ is the unique strong solution of \eqref{Eq:0.2} with $X(0,x) = x$ a.s.

The Markov process $(X(t,x),t\geq 0, x\in K\cap H_c)$ is continous and has $P^c$ for transition semigroup which is strong Feller on $H_c$.

For all $x \in K \cap H_c$ and $0=t_0<t_1<\cdots <t_m$, $(X(t_i,x) , i =1, \dots, n)$ is the limit in distribution of $(X^n(t_i,x))_{i =1, \dots, m}$.

Finally $\nu_c$ is an invariant measure for $P^c$.

\end{Th}

In all the article, C denotes a constant which may depend on $T$ and its value may change from one line to another.

\section{Proof of Theorem \ref{Th:Existence}}\label{S:2}

\subsection{Pathwise uniqueness}\label{SS:2.1}
We first prove that for any pair $(X^{i},\eta_{+}^{i},\eta_{-}^{i},W), i=1,2$, of weak solutions of \eqref{Eq:0.2} defined on the same probability space with the same driving noise $W$ and with $X^{1}(0) = X^{2}(0)$, we have $\left(X^{1},\eta_{+}^{1},\eta_{-}^{1}\right)=\left(X^{2},\eta_{+}^{2},\eta_{-}^{2}\right)$. This pathwise uniqueness will be used in the next subsection to construct stationary strong solutions of \eqref{Eq:0.2}.
\begin{Prop}
Let $x \in \mathcal{C}\left([0,1];[-1,1]\right)$. Let $\left(X^{i},\eta_{+}^{i},\eta_{-}^{i},W\right), i=1,2$ be two weak solutions of \eqref{Eq:0.2} with $X^{1}(0) = X^{2}(0)=x$. Then $\left(X^{1},\eta_{+}^{1},\eta_{-}^{1}\right)=\left(X^{2},\eta_{+}^{2},\eta_{-}^{2}\right)$.
\end{Prop}
\Proof We use the following Lemma from \cite{GOUDENEGE}.
\begin{Le}\label{Le:2.1}
Let $\zeta$ be a finite measure on $O_{\delta,T}$ and $V\in\mathcal{C}(O_{\delta,T})$. Suppose that there exists a positive continuous function $c_{T} : [0,T] \rightarrow \R^+$ such that :
\begin{enumerate}
\item[i)] for all $r\in[\delta,T]$, for all $h\in\mathcal{C}([0,1])$, such that $\bar{h}=0$, $\langle h,\zeta\rangle_{O_{r,T}}=0$,
\item[ii)] for all $r\in[\delta,T]$, $\overline{V(r,\cdot)} = c_{T}(r)$ with $\langle V,\zeta\rangle_{O_{r,T}}=0$,
\end{enumerate}
then $\zeta$ is the null measure.
\end{Le}

Let $Y(t) = X^1(t,x)-X^2(t,x)$, $\zeta_{+} = \eta_{+}^1-\eta_{+}^2$ and $\zeta_{-} = \eta_{-}^1-\eta_{-}^2$, $Y$ is the solution of the following equation:
\BE
\left\{\begin{array}{l}
\dd Y=-\frac{1}{2}A\left(AY+\left(f(X^1)-f(X^2)\right) + \zeta_{-} - \zeta_{+}\right)\ddt,\\
\\
Y(0)=0.
\end{array}\right.
\EE
Taking the scalar product in $H$ with $Y^N=Q_NY$ and integrating in time, 
we obtain since $Y$ has zero average:
\BE\label{Eq:2.2}
|Y^N(t)|^2_{-1}-|Y^N(\delta)|^2_{-1}=  -\int_\delta^t\left(|Y^N(s)|^2_{1}-\langle f(X^1)-f(X^2),Y^N\rangle\right)\dds+
\langle \zeta_{-}-\zeta_{+}, Y^N\rangle_{O_{\delta,t}}.
\EE
For all $s\in[\delta,t]$, 
$$
\begin{array}{l}
\langle Y^N(s) - Y(s), f(X^1(s,x))-f(X^2(s,x))\rangle\\
\\
\leq \|Y^N(s) - Y(s)\|_{L^\infty([0,1])}\|f(X^1(s,x))-f(X^2(s,x))\|_{L^1([0,1])},
\end{array}
$$
where $\|\cdot\|_{L^\infty([0,1])}$ and $\|\cdot\|_{L^1([0,1])}$ are the classical norms on the space $[0,1]$. The latter term converges to zero since $Y^N(s)$ converges uniformly to $Y(s)$ on $[0,1]$. Since $f(x)-\lambda x$ is nonincreasing,
\BEAS
\left(\langle Y(s), f(X^1(s,x))-f(X^2(s,x))\rangle\right) &=& \left(\langle Y(s), f(X^1(s,x))-f(X^2(s,x)) - \lambda Y(s)\rangle\right)\\
&&+ \left(\langle Y(s), \lambda Y(s)\rangle\right)\\
&\leq & \lambda |Y(s) |_{0}^2.
\EEAS
Taking the limit in \eqref{Eq:2.2} as $N$ grows to infinity, we obtain:
\BD
|Y(t)|_{-1}^2 - |Y(\delta)|_{-1}^2\leq \big\langle Y,\zeta_{-}-\zeta_{+}\big\rangle_{O_{\delta,t}} + \lambda \int_{\delta}^t |Y(s) |_{0}^2 \ \dds.
\ED
We now write
$$
\begin{array}{l}
 \big\langle Y,\zeta_{-}-\zeta_{+}\big\rangle_{O_{\delta,t}}\\
=  \big\langle 1+X^1,\eta_{-}^1\big\rangle_{O_{\delta,t}} -\big\langle 1+ X^2,\eta_{-}^1\big\rangle_{O_{\delta,t}} - \big\langle 1+X^1,\eta_{-}^2\big\rangle_{O_{\delta,t}} +\big\langle 1+X^2,\eta_{-}^2\big\rangle_{O_{\delta,t}}\\
+ \big\langle 1-X^1,\eta_{+}^1\big\rangle_{O_{\delta,t}} -\big\langle 1-X^2,\eta_{+}^1\big\rangle_{O_{\delta,t}} - \big\langle 1-X^1,\eta_{+}^2\big\rangle_{O_{\delta,t}} +\big\langle 1-X^2,\eta_{+}^2\big\rangle_{O_{\delta,t}}\\
\le 0
\end{array}
$$
by the contact condition and the positivity of the measures. 
It follows:
\BD
|Y(t)|_{-1}^2 - |Y(\delta)|_{-1}^2 \leq \lambda \int_{\delta}^t |Y(s) |_{0}^2 \ \dds.
\ED
By Gronwall Lemma, and letting $\delta \rightarrow 0$, we have $|Y(t)|_{-1}=0$ for all $t\geq0$. Since $\bar Y(t)=0$, we deduce $X^1(t,x) = X^2(t,x)$ for all $t\geq0$.
Moreover, with the definition of a weak solution, we see that :
\BD
\text{ for all } h\in D(A^2),\qquad  \big\langle Ah,\zeta_{+}-\zeta_{-}\big\rangle_{O_{\delta,t}} = 0.
\ED
By density, we obtain that $\zeta:=\zeta_{-}-\zeta_{+}$ and $V := (1-X^1)(1+X^1) = (1-X^2)(1+X^2)$ satisfy the hypothesis of Lemma \ref{Le:2.1}, and therefore $\zeta=\zeta_{-}-\zeta_{+}$ is the null measure. And since $\zeta_{-}$ and $\zeta_{+}$ have disjoint supports, then $\zeta_{-}$ and $\zeta_{+}$ are the null measure, i.e. $\eta_{-}^1=\eta_{-}^2$ and $\eta_{+}^1=\eta_{+}^2$.\EProof

\subsection{Convergence of invariants measures}\label{SS:2.2}
We know (see \cite{MR2349572}) that $\mu_c$ is the law of $Y^c= \B - \overline{\B} +c$, where $B$ is brownian motion. 
Then for $0\leq c < 1$, we remark the following inclusion : 
\BD
\{\B_{\theta} \in \left[\frac{c-1}{2},\frac{1-c}{2}\right], \text{ for all } \theta \in [0,1]\} \subset \{Y^c \in K\},
\ED
and we have a similar result for $-1< c \leq 0$. Therefore $\mu_c(K)>0$ with $-1<c<1$. Let us define $U$ the potential associated to the function $f$ :
\BD
U(x)=\left\{\begin{array}{rl}
\int_{0}^1F(x(\theta))\ddtheta &\text{ if } x \in K,\\
+\infty &\text{ else}.
\end{array}\right.
\ED
We have the following result : 
\begin{Prop} For $-1<c<1$,
\BD
\nu^n_{c} \rightharpoonup \nu_c := \frac{1}{Z_c}\exp^{-U(x)}\1_{x\in K}\mu_c(\ddx),\text{ when } n \rightarrow +\infty,
\ED
where $Z_{c}$ is a normalization constant.
\end{Prop}
\Proof
Let $\psi \in \mathcal{C}_{b}(L^2,\R)$. We want to prove that 
\BE\label{Eq:2.3}
\int_{H} \psi(x)\exp({-U^n(x)})\mu_{c}(\ddx) \mathop{\longrightarrow}_{n\rightarrow+\infty} \int_{H} \psi(x)\exp({-U(x)})\1_{x\in K}\mu_{c}(\ddx).
\EE
We first prove, 
\BE\label{Eq:2.4}
\exp({-U^n(x)}) \mathop{\longrightarrow}_{n\rightarrow+\infty} \exp({-U(x)})\1_{x\in K},\; \mu_c\; a.s.
\EE
Since $\mu_c(C([0,1]))=1$, we can restrict our attention to $x\in C([0,1])$. 
Then if $x\notin K$ there exists $\delta_{x}>0$ such that $m(\{\theta\in [0,1] : x(\theta)\leq-1-\delta_{x}\})>0$ or $m(\{\theta\in [0,1] : x(\theta)\geq1+\delta_{x}\})>0$, $m$ being the Lebesgue measure. Suppose $m(\{\theta\in [0,1] / x(\theta)\leq-1-\delta_{x}\})>0$, then we have since $\tilde F^n(x)=F^n(x)+\frac\lambda2 x^2$ is positive and non increasing on $(-\infty,-1)$
\BEAS
0\leq\exp({-U^n(x)}) &\leq& \exp\left(-\int_{0}^1\tilde F^n(x(\theta))
1_{\{x\leq -1-\delta_{x}\}}-\frac\lambda2 x(\theta)^2\ddtheta\right)\\
&\leq& \exp\left(-\int_{0}^1\tilde F^n(-1-\delta_{x})1_{\{x\leq -1-\delta_{x}\}}-\frac\lambda2 x(\theta)^2\ddtheta\right)\\
&\leq& \exp\left(-\tilde F^n(-1-\delta_{x})m(\{x\leq -1-\delta_{x}\})
+\int_{0}^1\frac\lambda2 x(\theta)^2\ddtheta\right).
\EEAS
And this latter term converges to zero as $n$ grows to infinity.\\
Now for $x \in K$, $F^n(x(\theta))$ converges to $F(x(\theta))$ almost everywhere as $n$  grows to infinity. Moreover $-\frac\lambda2 x(\theta)^2\leq F^n(x(\theta)) \leq \ln 2$, and by the dominated convergence Theorem, we deduce \eqref{Eq:2.4}.
Finally, \eqref{Eq:2.3} follows again by dominated convergence Theorem.\\

\subsection{Existence of stationary solutions}\label{SS:2.3}
In this section, we prove the existence of stationary solutions of equation \eqref{Eq:1.1} and that they are limits of stationary solutions of \eqref{Eq:1.2}, in some suitable sense. Fix $-1<c<1$ and consider the unique (in law) stationary solution of \eqref{Eq:1.2} denote $\hat{X}_{c}^{n}$ in $H_{c}$. We are going to prove that the laws of $\hat{X}_{c}^{n}$ weakly converge as $n$ grows to infinity to a stationary strong solution of \eqref{Eq:0.2}.\\
\begin{Prop}\label{Prop:2.3}
Let $-1<c<1$ and $T>0$, $\hat{X}^n_{c}$ converges in probability as $n$ grows to infinity to a process $\hat{X}_{c}$ in $\mathcal{C}(O_{T})$. Moreover $f(\hat{X}_{c}) \in L^1(O_T)$ almost surely, and setting
\BD
\dd\eta^n_{+} = -f^n(\hat{X}^n_{c}(t,\theta))\1_{\hat{X}^n_{c}(t,\theta)> 0}\ddt\ddtheta+f(\hat{X}_{c}(t,\theta))\1_{0<\hat{X}_{c}(t,\theta)\leq 1}\ddt\ddtheta,
\ED
and
\BD
\dd\eta^n_{-} = f^n(\hat{X}^n_{c}(t,\theta))\1_{\hat{X}^n_{c}(t,\theta)\leq 0}\ddt\ddtheta-f(\hat{X}_{c}(t,\theta))\1_{-1\leq\hat{X}_{c}(t,\theta)\leq 0}\ddt\ddtheta,
\ED
then $(\eta^n_{+},\eta^n_{-})$ converges in probability to $(\eta_{+},\eta_{-})$ such that $(\hat{X}_{c},\eta_{+},\eta_{-},W)$ is a stationary strong solution of \eqref{Eq:0.2}.
\end{Prop}
\Proof
Proceeding exactly as in \cite{MR2349572} (see Lemma 5.2), we prove that the laws of $(\hat X^n_c,W^n)_{n\in\N}$ are tight in $C(O_T)\times C([0,T];V_\gamma)$, $\gamma < -1/2$. We have set $W^n=W$, $n\in\N$. We therefore can extract convergent subsequences. Let $(\hat X^{n_k}_c, W^{n_k})_{k\in\N}$ be 
such a subsequence. Using Skohorod theorem, one may find a probability space and a sequence of random variables $(\tilde X^{k}_c, \mathcal{W}^k)_{k\in\N}$ on this probability space with the same laws as $(\hat X^{n_k}_c,W^{n_k})_{k\in\N}$ which converge almost
surely. 

Below, we  show in Step 1 that its limit $\tilde X_c$ satisfies $f(\tilde{X}_{c}) \in L^1(O_{T})$ almost surely. Then in Step 2, we prove that  the measures $\tilde\eta^{k}_{\pm}$, defined as above with $\hat X^{n_k}_c$ replaced by $\tilde X^k_c$, converges to two positive measures $\tilde\eta_{\pm}$ and that $(\tilde{X}_{c},\tilde\eta_{+},\tilde\eta_{-})$ is a weak solution in the probabilistic sense. It then remains to use pathwise uniqueness to conclude in Step 3. In this proof, we only treat the case $\lambda=0$. This assumption is not essential at all but lightens the computations. For $\lambda \ne 0$, an extra term has to be taken into account. It is very easy to deal with. \\

\Step[1] \\
Applying Ito formula to $|Q_N\hat X_c^n(t)|_{-1}^2$, we obtain
$$
\begin{array}{l}
\ds |Q_N\hat X_c^n(T)|_{-1}^2-|Q_N\hat X_c^n(0)|_{-1}^2+Ê\int_0^T|Q_N\hat X_c^n(t)|_{1}^2\ddt 
-2\int_{O_T} f_n(\hat X_c^n)\left(Q_N\hat X_c^n -c\right) \dds\ddtheta\\
\ds = 2 \int_0^T(Q_N\hat X_n^c,B\dd W(s)) +   T\, \Tr(Q_NB) 
\end{array}
$$
Note that 
$$
\E\left( \left(\int_0^T (Q_N\hat X_c^n,B\dd W(s))\right)^2 \right)= \E\int_0^T |Q_N\hat X_c^n|^2_{-1}\dds 
=T\int_H |Q_Nx|^2_{-1} \nu_c^n(\ddx) \le C\, T
$$
We set 
$$
\begin{array}{rl}
\varphi_{n}^N=  |Q_N\hat X_c^n(T)|_{-1}^2-|Q_N\hat X_c^n(0)|_{-1}^2&+Ê\int_0^T|Q_N\hat X_c^n(t)|_{1}^2\ddt\\ 
\\
&-2\int_{O_T} f_n(\hat X_c^n)\left(Q_N\hat X_c^n -c\right) \dds\ddtheta- T\, \Tr(Q_NB) 
\end{array}
$$
and deduce
$$
\P( |\varphi_n^N|\ge M)\le \frac{C\, T}{M^2} .
$$
Thus, for all $N\inÊ\N$, the laws of $(\varphi_n^N)_{n\in\N}$ are tight. Therefore the laws of 
$(\hat X^n_c,W^n,(\varphi_n^N)_{N\in\N})_{n\in\N} $ are tight and using 
Skohorod theorem on this sequence, we can assume that $\tilde X^{k}_c$, $\mathcal{W}^k$ and, for $N\in\N$, $\tilde \varphi_{k}^N$ converge
almost surely. We have defined $\tilde \varphi_{k}^N$ as above with $\tilde X^{k}_c$ instead of $\hat X^n_c$. 
In particular, $\tilde \varphi_{k}^N$ is bounded almost surely:
$$
\begin{array}{rl}
\ds |Q_N\tilde X^k_{c}(T)|_{-1}^2-|Q_N\tilde X^k_{c}(0)|_{-1}^2&+2\int_0^T|Q_N\tilde X^k_{c}(t)|_{1}^2\ddt 
\\
\\&-2\int_{O_T} f_{k}(\tilde X^k_{c})\left(Q_N\tilde X^k_{c} -c\right) \dds\ddtheta- T\, \Tr(Q_NB)\\
\\
&\ds \le C(N,T,c)
\end{array}
$$
where $C(N,T,c)$ is random. The first three terms are clearly also bounded almost surely. This uses 
the fact that $Q_N$ is a bounded operator from $H$ to $V_1$. Since $Q_N$ has finite dimensional range,
we obtain
\BE\label{Eq:2.5}
-\int_{O_T} f_{n_k}(\tilde X_c^{k})\left(Q_N\tilde X_c^{k} -c\right) \dds\ddtheta
 \le C(N,T,c)
\EE
for a different random constant $C(N,T,c)$.

Let us choose $\epsilon_0 =\min\left\{\frac{1-c}{4},\frac{1+c}4\right\}$ and take $N\in\N$ such that
$$
|Q_N\tilde X_c -\tilde X_c|_{C(O_T)}\le \frac12 \epsilon_0 
$$
and $K_0$ such that for $k\ge K_0$ 
$$
|\tilde X_c^{k}-\tilde X_c|_{C(O_T)}\le \frac 14\epsilon_0 .
$$
Then, for all $k\ge K_0$,
$$
|Q_N\tilde X_c^{k} -\tilde X_c^{k}|_{C(O_T)}\le  \epsilon_0 .
$$
Moreover, if $\tilde X_c^{k}\ge \frac{1+c}2$ then $f_{n_k}(\tilde X_c^{k})\le 0$ and
$$
Q_N\tilde X_c^{k}-c\ge -\epsilon_0 + \frac{1+c}2 -c
\ge \frac{1-c}4\ge  \epsilon_0.
$$
 Similarly, if $\tilde X_c^{k}\le \frac{-1+c}2$ then $f_{n_k}(\tilde X_c^{k})\ge 0$ and
$$
Q_N\tilde X_c^{k}-c \le -  \epsilon_0.
$$
Finally, noticing that $f_n$ is uniformly bounded by a constant $K(c)$ on $[\frac{-1+c}2, \frac{1+c}2]$,
we deduce
\BEAS
\int_{O_T}|f_{n_k}(\tilde X_c^{k})|\dds\ddtheta & \le& - \frac1{\epsilon_0}\int_{\tilde X_c^{k}\ge \frac{1+c}2} 
f_{n_k}(\tilde X_c^{k})\left(Q_N\tilde X_c^{k} -c\right) \dds\ddtheta \\
\\
&& -\frac1{\epsilon_0}\int_{\tilde X_c^{k}\le \frac{-1+c}2} 
f_{n_k}(\tilde X_c^{k})\left(Q_N\tilde X_c^{k} -c\right) \dds\ddtheta
+ K(c)\\
\\
&\leq&  -\frac1{\epsilon_0}\int_{O_T} 
f_{n_k}(\tilde X_c^{k})\left(Q_N\tilde X_c^{k} -c\right) \dds\ddtheta\\
\\
&& +\frac18(\max\{1-c,1+c\})^2K(c) + K(c).
\EEAS
Thanks to \eqref{Eq:2.5}, we obtain 
\BE\label{Eq:2.6}
\int_{O_T}|f_{n_k}(\tilde X_c^{k})|\dds\ddtheta \le C(N,T,c),
\EE
where the value of the random constant $C(N,T,c)$ has again changed. It easily deduced from this 
uniform bound that $|\tilde X_c|\le 1$ almost everywhere with respect to $t$ and $\omega$ and by 
Fatou Lemma that $f(\tilde X_c)\in L^1(O_T)$ almost surely.
\EStep[1]

 \Step[2] \\

Let now 
$\xi^k$ be the following measure on $O_{T}$: 
\BD
\dd\xi^k:=-f^{n_{k}}(\tilde X_c^{k}(t,\theta))\ddt\ddtheta.
\ED
and $\xi_{+}^k$ and $\xi_{-}^k$ the positive and negative parts: 
\BD
\dd\xi_{+}^k:=-f^{n_{k}}(\tilde X_c^{k}(t,\theta))\1_{\tilde X_c^{k} > 0}\ddt\ddtheta,\;
\dd\xi_{-}^k:=f^{n_{k}}(\tilde X_c^{k}(t,\theta))\1_{\tilde X_c^{k}\leq 0}\ddt\ddtheta.
\ED
By step 1, $f(\tilde X_c) \in L^1(O_{T})$ and we can define the following measure:
\BD
\dd\lambda := -f((\tilde X_c(t,\theta))\1_{-1\leq \tilde X_c\leq 1}\ddt\ddtheta,
\ED
and the positive and negative parts:
\BD
\dd\lambda_{+} := -f((\tilde X_c(t,\theta))\1_{0< \tilde X_c\leq 1}\ddt\ddtheta,\quad \dd\lambda_{-} := f((\tilde X_c(t,\theta))\1_{-1\leq \tilde X_c\leq 0}\ddt\ddtheta.
\ED
By \eqref{Eq:2.6}, $f^{n_{k}}(\tilde X_c^{k})-f(\tilde X_c)$ is bounded in $L^1(O_T)$. We deduce that $\xi^k$ 
has a subsequence $\xi^{k_\ell}$which converges to a measure $ \zeta$. 
Note that this subsequence may depend on the random parameter $\omega$. We set $\tilde\eta= \zeta -\lambda$. 

Thanks to Fatou Lemma we have the following inequality for all $h \in \mathcal{C}(O_{T})$ nonnegative:
\BEAS
\int_{O_{T}} h(s,\theta) \big[-f(\tilde X_c(s,\theta))\1_{0< \tilde X_c\leq 1} \big]\dd s \dd\theta &\!\!\!=\!\!\!& \int_{O_{T}}\mathop{\liminf}_{\ell\rightarrow+\infty} \big[-h(s,\theta)f^{n_{{k_\ell}}}(\tilde X_c^{{k_\ell}}(s,\theta))\1_{0< \tilde X_c^{{k_\ell}} \leq 1}\big]\dds \ddtheta\\
&\!\!\!\leq\!\!\!&\mathop{\liminf}_{\ell\rightarrow+\infty} \int_{O_{T}}\big[-h(s,\theta)f^{n_{{k_\ell}}}(\tilde X_c^{{k_\ell}}(s,\theta))\1_{0< \tilde X_c^{{k_\ell}}\leq 1}\big]\dds \ddtheta.
\EEAS
Therefore $\tilde\eta^{{k_\ell}}_+=\xi_{+}^{k_\ell}-\lambda_{+}$ converges to a positive measure. 
Similarly, $\tilde\eta^{{k_\ell}}_-=\xi_-^{k_\ell}-\lambda_-$ converges to a positive measure. 
It follows:
\BD
\xi_{+}^{k_\ell}-\lambda_{+} \rightharpoonup \tilde\eta_+\quad \text{ and }\quad\xi_{-}^{k_\ell}-\lambda_{-} \rightharpoonup \tilde\eta_-,
\ED
where $\tilde\eta_+$ and $\tilde\eta_-$ are the positive and negative parts of $\tilde\eta$.

Let us now show that the contact conditions holds for $\left(1-\tilde X_c,\tilde\eta_+\right)$ 
and $\left(1+\tilde X_c,\tilde\eta_-\right)$.
Let us define the following measures for $\varepsilon>0$ and  $k\in\N$.
\BD
\begin{array}{ll}
\dd\xi^{k}_{+,\varepsilon} := -f^{n_{k}}(\tilde X_c^{k}(t,\theta))\1_{1-\varepsilon\leq \tilde X_c^{k}}\ddt\ddtheta,
&\dd\tau^{k}_{+,\varepsilon} := -f^{n_{k}}(\tilde X_c^{k}(t,\theta))\1_{0< \tilde X_c^{k}<1-\varepsilon}\ddt\ddtheta,\\
\\	
\dd\lambda_{+,\varepsilon} := -f(\tilde X_c(t,\theta))\1_{1-\varepsilon\leq \tilde X_c}\ddt\ddtheta, &\dd\tau_{+,\varepsilon} := -f(\tilde X_c(t,\theta))\1_{0< \tilde X_c<1-\varepsilon}\ddt\ddtheta.
\end{array}
\ED
Clearly $\tau_{+,\varepsilon}^k$ converges to $\tau_{+,\varepsilon}$, it follows
\BEAS
\limsup_{\ell\rightarrow +\infty}\left\langle 1-\tilde X_c^{{k_\ell}}, \xi_{+}^{k_\ell}-\lambda_{+}\right\rangle_{O_{T}}&=&\limsup_{\ell\rightarrow +\infty}\left( \big\langle 1-\tilde X_c^{{k_\ell}},\xi^{{k_\ell}}_{+,\varepsilon}\big\rangle_{O_{T}}\!\!-\big\langle 1-\tilde X_c^{{k_\ell}},\lambda_{+,\varepsilon}\big\rangle_{O_{T}}\right.\\
&&+\left.  \big\langle 1-\tilde X_c^{{k_\ell}},\tau^{k_\ell}_{+,\varepsilon}\big\rangle_{O_{T}}\!\!-\big\langle 1-\tilde X_c^{{k_\ell}},\tau_{+,\varepsilon}\big\rangle_{O_{T}}\right)\\
\\
&=&\limsup_{\ell\rightarrow +\infty}\left(\int_{O_{T}} \!\!\left(\tilde X_c^{{k_\ell}}-1\right)f^{n_{{k_\ell}}}(\tilde X_c^{{k_\ell}})\1_{1-\varepsilon\leq \tilde X_c^{{k_\ell}}}\ddt\ddtheta\right.\\
&&\left. + \int_{O_{T}} \!\!\left(1-\tilde X_c^{{k_\ell}}\right)f(\tilde X_c)\1_{1-\varepsilon\leq \tilde X_c}\ddt\ddtheta\right)\\
\\
&\leq&\limsup_{\ell\rightarrow +\infty}	\left(\int_{O_{T}} \left(\tilde X_c^{{k_\ell}}-1\right)f^{n_{{k_\ell}}}(\tilde X_c^{{k_\ell}})\1_{1-\varepsilon\leq \tilde X_c^{{k_\ell}}\leq 1}\ddt\ddtheta\right)\\
&& + \limsup_{\ell\rightarrow +\infty}	\left(\int_{O_{T}} \left(1-\tilde X_c^{{k_\ell}}\right)^{-}f(\tilde X_c)\1_{1-\varepsilon\leq \tilde X_c}\ddt\ddtheta \right)
\EEAS
Since $(1-\tilde X_c^{{k_\ell}})^-$ converges uniformly to zero, we deduce:
\BEAS
\limsup_{\ell\rightarrow +\infty} \ \langle 1-\tilde X_c^{{k_\ell}},  \xi_{+}^{k_\ell}-\lambda_{+}\big\rangle_{O_{T}}&\leq&T \sup_{x\in[1-\varepsilon,1]}\left|(x-1)f(x)\right|\\
\\
&\leq&-T\varepsilon\ln\left(\frac{\varepsilon}{2-\varepsilon}\right).
\EEAS
Letting $\varepsilon \to 0$, we obtain the first contact condition since the left hand side clearly converges to $\langle1-\tilde X_c,\tilde\eta_+\rangle$. The second is obtained similarly.

We now prove that $\xi^k-\lambda$ does not have more than one limit point so that in fact the whole sequence converge to $\tilde\eta$. Let $\tilde\eta_i$, $i=1,2$ be two limit points.

For all $h \in D(A^2)$ and for all $0\leq t \leq T$:
\BEAS
\big\langle Ah, \xi^k- \lambda\big\rangle_{O_{t}} &=&
\langle \tilde X_c^{k}(t,.),h\rangle -\langle x,h\rangle + \int_{O_{t}} \tilde X_c^{k}(s,\theta)A^2h(\theta) \dds\ddtheta\\
&& + \int_{0}^{t}\langle Bh, \dd\mathcal{W}^k\rangle + \int_{O_{t}} f(\tilde X_c(s,\theta))Ah(\theta) \dds\ddtheta.
\EEAS
We deduce 
\BEAS
\big\langle Ah,  \tilde\eta_i\big\rangle_{O_{t}}&=&
-\langle \tilde X_c(t,.),h\rangle +\langle x,h\rangle - \int_{O_{t}} \tilde X_c(s,\theta)A^2h(\theta) \dds\ddtheta\\
&& - \int_{0}^{t}\langle Bh, \dd\mathcal{W}\rangle - \int_{O_{t}} f(\tilde X_c(s,\theta))Ah(\theta) \dds\ddtheta.
\EEAS
And by density
$$ 
\big\langle h,  \tilde\eta_1\big\rangle_{O_{t}}=\big\langle h,  \tilde\eta_2\big\rangle_{O_{t}}
$$
for any $h\in C([0,1])$ such that $\bar h=0$. Since by the contact condition
$$
\big\langle (1-\tilde X_c)(1+\tilde X_c),  \tilde\eta_1\big\rangle_{O_{t}}=\big\langle (1-\tilde X_c)(1+\tilde X_c),  \tilde\eta_2\big\rangle_{O_{t}}.
$$
We deduce from Lemma \ref{Le:2.1} that $\tilde\eta_1=\tilde\eta_2$.
\EStep[2]

\indent\Step[3]\\
We use a result form \cite{MR1392450} that allows to get the convergence of the approximated solutions in probability in any space in which these approximated solutions are tight.
\begin{Le}\label{Le:2.2}
Let $\{Z_{n}\}_{n\geq1}$ be a sequence of random elements on a Polish space $E$ endowed by its borel $\sigma$-algebra. Then $\{Z_{n}\}_{n\geq1}$ converges in probability to an $E$-valued 
random element if and any if from every pair of subsequences $\{(Z_{n_{k}^1},Z_{n_{k}^2})_{k\geq1}$, one 
can extract a subsequence which converges weakly to a random element supported on the diagonal 
$\{(x, y) \in E \times E, x = y\}$.
\end{Le}
Assume $(n^1_{k})_{k\in\N}$ and $(n^1_{k})_{k\in\N}$ are two arbitrary subsequences. Clearly, the process $\left(\hat{X}^{n^1_{k}}_{c},\hat{X}^{n^2_{k}}_{c},W^k\right)$ is tight in a suitable space. By Skorohod's theorem, we can find a probability space and a sequence of processes $\left(\tilde{X}^{1,k}_{c},\tilde{X}^{2,k}_{c},\mathcal{W}^k\right)$ such that $\left(\tilde{X}^{1,k}_{c},\tilde{X}^{2,k}_{c},\mathcal{W}^k\right) \rightarrow \left(\tilde{X}^{1}_{c},\tilde{X}^{2}_{c},\mathcal{W}\right)$ almost surely  and 
$\left(\tilde{X}^{1,k}_{c},\tilde{X}^{2,k}_{c},\mathcal{W}^k\right)$ has the same distribution as $\left(\hat{X}^{n^1_{k}}_{c},\hat{X}^{n^2_{k}}_{c},W^k\right)$ for all $k\in\N$. In the Skorohod's space, the approximated measures respectively converge to two contact measures $\tilde\eta_{1}$ and $\tilde\eta_{2}$. By the second step, $(\tilde{X}^{1}_{c},\tilde\eta_{1},\mathcal{W})$ and $(\tilde{X}^{2}_{c},\tilde\eta_{2},\mathcal{W})$ are both weak solutions of \eqref{Eq:0.2}. By uniqueness, necessarily $\tilde{X}^{1}_{c}=\hat{X}^{2}_{c}$ and $\tilde\eta_{1}=\tilde\eta_{2}$. Therefore the subsequence $\left(\hat{X}^{n^1_{k}}_{c},\hat{X}^{n^2_{k}}_{c}\right)_{k\in\N}$ converges in distribution to a process supported on the diagonal. We use Lemma \ref{Le:2.2} to prove that the sequence $(\hat{X}^{n}_{c})$ converges in 
probability to a process $\hat{X}_{c}$. Clearly  $\hat{X}_{c}$ is stationary. Reproducing the argument of Step 1 and Step 2, we prove that it is a strong solution of \eqref{Eq:0.2} and the convergence of the contact measures. \EStep[3]

\subsection{Convergence of the semigroup}\label{SS:2.4}
First we state the following result which is a corollary of Proposition \ref{Prop:2.3}.
\begin{Coro}\label{Co:2.1}Let $c>0$.
\begin{enumerate}
\item[i)] There exists a continuous process $(X(t,x),t\geq 0,x\in K\cap H_c)$ with $X(0,x) = x$ and a set $K_0$ dense in $K\cap H_c$, such that for all $x \in K_0$ there exists a unique strong solution of equation \eqref{Eq:0.2} given by $\left(\left(X(t,x)\right)_{t\geq0},\eta^x_{+},\eta^x_{-},W\right)$.
\item[ii)] The law of $\left(X(t,x)_{t\geq 0},\eta^x_{+},\eta^x_{-}\right)$ is a regular conditional distribution of the law of $\left(\hat{X}_c,\eta_{+},\eta_{-}\right)$ given $\hat{X}_{c}(0)=x \in K \cap H_c$.
\end{enumerate}
\end{Coro}
\Proof By Proposition \ref{Prop:2.3}, we have a stationary strong solution $\hat{X}_{c}$ in $H_c$, such that $W$ and $\hat{X}_{c}(0)$ are independent. Conditioning $\left(\hat{X}_{c},\eta_{+},\eta_{-}\right)$ on the value of $\hat{X}_{c}(0)=x$, with $c=\overline{x}$, we obtain for $\nu_c$-almost every $x$ a strong solution that we denote $\left(X(t,x),\eta^x_{+},\eta^x_{-}\right)$ for all $t\geq0$ and for all $x\in K\cap H_{c}$. This process is the desired process. Indeed, since the support of $\nu_c$ is $K\cap H_c$, we have a strong solution for a dense set $K_0$ in $K \cap H_c$.\\
Notice that all processes $\left(X(t,x)\right)_{t\geq 0}$ with $x \in K_0$ are driven by the same noise $W$ and are continuous with values in $H$. Moreover, we have the following obvious identity:
\BD
|X^n(t,x)-X^n(t,y)|_{-1} \leq e^{\lambda t} |x-y|_{-1},\; x,y \in L^2_c,\; t\ge 0,
\ED
and by density we obtain a continuous process $\left(X(t,x)\right)_{t\geq 0}$ in $H_c$ for all $x \in K \cap H_c$.
\EProof
\begin{Prop}\label{Prop:2.4}
Let $c>0$, for all $\phi \in \mathcal{C}_b(H)$ and $ x \in K \cap H_c$:
\BE\label{Eq:2.7}
\lim_{n \rightarrow +\infty} P^{n,c}_t\phi(x) = \E[\phi(X(t,x))]=:P^c_t\phi(x).
\EE
Moreover the Markov process $(X(t,x),t\geq 0, x\in K\cap H_c)$ is strong Feller and its transition semigroup $P^c$ is such that:
\BE\label{Eq:2.8}
|P^c_t\phi(x)-P^c_t\phi(y)|\leq \frac{2e^{\lambda^2t/4}}{\lambda\sqrt{t}}|x-y|_{-1},\quad \text{ for all }x,y \in K \cap H_c,\text{for all } t>0.
\EE
\end{Prop}
\Proof 
Since $(\nu_{c}^n)_{n\geq 1}$ is tight in $H_{c}$, then there exists an increasive sequence of compact sets $(J^p)_{p\in\N}$ in $H$ such that:
\BD
\lim_{p\rightarrow +\infty} \sup_{n\geq1} \nu_{c}^n(H\setminus J^p) = 0.
\ED
Set $J:= \mathop{\cup}_{p\in\N} J^p\cap K$. Since the support of $\nu_{c}$ is in $K\cap H_{c}$ and $\nu_{c}(J)=1$, then $J$ is dense in $K\cap H_{c}$.
Fix $t>0$, by \eqref{Eq:1.3}, for any $\phi \in \mathcal{C}_b(H)$ :
\BD
\sup_{n \in \N}(\|P^{n,c}_t\phi\|_{\infty} + [P^{n,c}_t\phi]_{Lip(H_c)}) < +\infty.
\ED
Let $(n_j)_{j\in\N}$ be any sequence in $\N$. With a diagonal procedure, by Arzel\`a-Ascoli Theorem, there exists $(n_{j_{l}})_{l\in\N}$ a subsequence and a function $\Theta_{t}:J\rightarrow \R$ such that:
\BD
\lim_{l\rightarrow +\infty} \sup_{x\in J^p}|P^{n_{j_{l}},c}_t\phi(x) - \Theta_{t}(x)|=0, \quad \text{ for all }p\in\N.
\ED
By density, $\Theta_{t}$ can be extended uniquely to a bounded Lipschitz function $\tilde{\Theta}_{t}$ on $K\cap H_{c}$ such that
\BD
\tilde{\Theta}_{t}(x)=\lim_{l\rightarrow + \infty}P^{n_{j_l},c}_t\phi(x),\quad\text{ for all } x \in K\cap H_c.
\ED
Note that the subsequence depends on $t$. Therefore, we have to prove that the limit defines a semigroup and does not depend on the chosen subsequence.\\
By Proposition \ref{Prop:2.3}, we have for all $\phi,\psi \in \mathcal{C}_b(H)$ :
\BEAS
\E\left[\psi\left(\hat{X}_{c}(0)\right)\phi\left(\hat{X}_{c}(t)\right)\right] &=&\lim_{l\rightarrow +\infty}\E\left[\psi\left(\hat{X}^{n_{j_{l}}}_{c}(0)\right)\phi\left(\hat{X}^{n_{j_{l}}}_{c}(t)\right)\right]\\
&=& \lim_{l\rightarrow +\infty}\int_H\psi(y)\E\left[\phi\left(\hat{X}^{n_{j_{l}}}_{c}(t)\right)\Big|\hat{X}^{n_{j_{l}}}_{c}(0)=y\right]\nu^{n_{j_{l}}}_c(\ddy)\\
&=& \lim_{l\rightarrow +\infty}\int_H\psi(y)P^{n_{j_{l}},c}_{t}\phi(y)\nu^{n_{j_{l}}}_c(\ddy)\\
&=& \int_H\psi(y)\tilde{\Theta}_{t}(y)\nu_c(\ddy).
\EEAS
Thus, by Corollary \ref{Co:2.1}, we have the following equality:
\BE\label{Eq:2.9}
\E\left[\phi\left(X(t,x)\right)\right]=\tilde{\Theta}_{t}(x), \quad \text{ for }\nu_c \text{-almost every } x.
\EE
Since $\E[\phi(X(t,.))]$ and $\tilde{\Theta}_{t}$ are continuous on $K\cap H_c$, and $\nu_c(K\cap H_c)=1$, the equality \eqref{Eq:2.9} is true for all $x \in K\cap H_c$. Moreover the limit does not depend on the chosen subsequence, and we obtain \eqref{Eq:2.7}. Letting $n\to\infty$ in \eqref{Eq:1.3}, we deduce \eqref{Eq:2.8}.
\EProof

\subsection{End of the proof of Theorem \ref{Th:Existence}}\label{SS:2.5}
We have proved that there exists a continous process $X$ which is a strong solution of equation \eqref{Eq:0.2}  for an $x$ in a dense space. In this section, we prove existence for an initial condition in $K\cap H_c$ with $c>0$. 

By Corollary \ref{Co:2.1} we have a process $(X(t,x),t\geq0,x\in K \cap H_c)$, such that for all $x$ in a set $K_0$ dense in $K \cap H_c$ we have a strong solution $\left(\left(X(t,x)\right)_{t\geq0},\eta^x_{+},\eta^x_{-},W\right)$ of \eqref{Eq:0.2} with initial condition $x$. By Proposition \ref{Prop:2.3}, the Markov process $X$ has transition semigroup $P^c$ on $H_c$.\\
The strong Feller property of $P^c$ implies that for all $x\in K \cap H_c$ and $s>0$ the law of $X(s,x)$ is absolutely continous with respect to the invariant measure $\nu_c$. Indeed, if $\nu_{c}(\Gamma)=0$, then $\nu_c(P_s^c(\1_{\Gamma})) = \nu_c(\Gamma)=0$. So $P_s^c(\1_{\Gamma})(x)=0$ for $\nu_c$-almost every $x$ and by continuity for all $x \in K \cap H_c$.\\
Therefore almost surely $X(s,x) \in K_0$ for all $s>0$ and $x \in K \cap H_c$. Fix $s>0$, denote for all $\theta \in [0,1]$:
\BD
\tilde{X} := t \mapsto X(t+s,x),
\tilde{W}(\cdot,\theta) := t \mapsto W(t+s,\theta)-W(s,\theta)),
\ED
and the measures $\tilde{\eta}^x_{\pm}$ such that for all $T>0$, and for all $h \in \mathcal{C}(O_{T})$:
\BD
\big\langle h, \tilde{\eta_{\pm}}^x\big\rangle_{O_{T}} := \int_{O_{s}^{T+s}} h(t-s,\theta)\eta^x_{\pm}(\ddt,\ddtheta).
\ED
So we have a process $\tilde{X} \in \mathcal{C}([0,T];H) \cap \mathcal{C}(O_{T})$ and two measures $\tilde{\eta}^x_{+}$ and $\tilde{\eta}^x_{-}$ on $O_{T}$ which is finite on $[\delta,T]\times[0,1]$ for all $\delta \geq0$, such that $\left((\tilde{X}(t,x))_{t\geq0},\tilde{\eta}^x_{+},\tilde{\eta}^x_{-},\tilde{W}\right)$ is a strong solution of \eqref{Eq:0.2} with initial condition $X(s,x)$. By continuity $X(s,x)\rightarrow x$ in $H$ as $s\rightarrow0$, so $\left((X(t,x))_{t\geq0},\eta^x_{+},\eta^x_{-},W\right)$ is a strong solution of \eqref{Eq:0.2} with initial condition $x$ in the sense of the definition \ref{De:1}.\\

\section{Ergodicity and mixing}

When $\lambda$ is small, it can be easily shown that $\nu_c$ is the unique invariant measure
and is ergodic. We now prove that this is in fact true for any $\lambda$. Note that since
$(P^c_t)_{t\ge 0}$ is Strong Feller, the results follows from Doob theorem if we prove that 
$(P^c_t)_{t\ge 0}$ is irreducible (see for instance \cite{MR1417491}). For additive noise driven SPDEs, this is often proved
by a control argument and continuity with respect to the noise. This latter property is not completely trivial
in our situation but we are able to adapt the argument.

\begin{Prop}
\label{p3.1}
For any $c\in (-1,1)$, the semigroup $(P^c_t)_{t\ge 0}$ is irreducible.
\end{Prop}

\Proof

Let $x,\, y\in C^\infty([0,1])$ be such that $|x|_{L^\infty(0,1)}\le 1-\delta$ and  $|y|_{L^\infty(0,1)}\le 1-\delta$
for some $\delta >0$ and $\bar x=\bar y=c$. We set 
$$
u(t)=\frac{t}{T} y +\left(1-\frac{t}{T}\right)x
$$
and define $g_0$ by
$$
g_0(\theta,t)=\int_0^\theta \left(\frac{1}{T}(y-x)+\frac12A(Au+f(u))\right)(\vartheta,t) \dd\vartheta 
$$
Then $g_0$ is in $C^{\infty}([0,T]\times [0,1])$, $g_0(t)\in D(B)$, $t\in [0,T]$,  and:
$$
\frac{d}{dt}u=-\frac12A(Au+f(u)) +Bg_0.
$$
Moreover
\BE\label{Eq:ergo1}
\frac{d}{dt}u=-\frac12A(Au+f_\delta(u)) +Bg_0
\EE
where $f_\delta$ is any Lipschitz function equal to $f$ on $[-1+\delta/2,1-\delta/2]$.

Let $X^\delta(\cdot,x)$ be the solution of \eqref{Eq:0.2} with $f$ replaced by $f_\delta$ and set
$Y^\delta(\cdot,x)=X^\delta(\cdot,x) -Z$, where $Z=Z(\cdot,0)$ is the solution of the linear equation
with $0$ as initial data. Then
$$
\frac{d}{dt}Y^\delta = -\frac{1}{2}A\left(AY^\delta +f_\delta(Y^\delta+Z)\right),\; Y^\delta(0,x)=x.
$$
Let also 
$$
z_0(t)=\int_0^t e^{-A^2(t-s)/2}Bg_0(s)\dds.
$$
Since the gaussian process $Z$ is almost surely continuous and has a non degenerate covariance, we clearly have
$$
\P\left(|Z-z_0|_{C(O_{T})} \le \varepsilon \right) >0
$$
for any $\varepsilon >0$. Let us denote by $Y^z$ the solution of 
\BE
\label{e3.2}
\frac{d}{dt}Y^z= -\frac{1}{2}A\left(AY^z +f_\delta(Y^z+z)\right),\; Y^z(0,x)=x.
\EE
We prove below that the mapping
$$
\Phi_\delta\;:\; z\mapsto Y^z
$$
is continuous from $C(O_{T})$ into $C(O_{T})$. Since $u=\Phi_\delta(z_0)+z_0$
and $X^\delta=\Phi_\delta(Z)+Z$, we
deduce that  there exists $\varepsilon$ such that
$$
\P\left(|X^\delta-u|_{C(O_{T})} \le \delta/2 \right) \ge \P\left(|Z-z_0|_{C(O_{T})} \le \varepsilon \right) >0
$$
Let us now observe that $|X^\delta-u|_{C(O_{T})} \le \delta/2 $ implies $|X^\delta|_{C(O_{T})} \le 1-\delta/2 $ so that $f_\delta(X^\delta) = f(X^\delta)$ and $X^\delta$ is fact solution of \eqref{Eq:0.2}. By pathwise uniqueness,
we deduce that $|X^\delta-u|_{C(O_{T})} \le \delta/2 $ implies $X^\delta=X$. It follows
$$
\P\left(|X-u|_{C(O_{T})} \le \delta/2 \right)\ge \P\left(|X^\delta-u|_{C(O_{T})} \le \delta/2 \right)>0
$$
In particular
$$
\P\left(|X(T,x)-y| \le \delta/2 \right)>0.
$$

If we assume now that $x,y\in H_c$ , we choose $\tilde x, \tilde y\in C^\infty([0,1])$ such that
$$
|x-\tilde x|\le \delta,\; |y-\tilde y|\le \delta, \; |\tilde x|_{L^\infty(0,1)}\le 1-\delta \mbox{ and }|\tilde y|_{L^\infty(0,1)}\le 1-\delta,
$$
and $\bar {\tilde x}=\bar{\tilde y}=c$.
We have 
$$
|X(T,x)-X(T,\tilde x)|\le e^{\lambda T}|x-\tilde x|.
$$
Therefore
$$ 
\P\left(|X(T,x)-y| \le \delta/2 + (1+e^{\lambda T})\delta \right) \ge \P\left(|X(T,\tilde x)-\tilde y| \le \delta/2 \right) >0.
$$
This proves the results. 

It remains to prove that $\Phi_\delta$ is continuous. This follows form the mild form 
of equation \eqref{e3.2}:
$$
Y^z(t)= e^{-tA^2/2}x +\int_0^t e^{-(t-s)A^2/2}Af_\delta (Y^z(s)+z(s))ds.
$$
It is classical that, for $t>0$, $Ae^{-tA^2/2}$ maps $C([0,1])$ into itself and
$$
\left| Ae^{-tA^2/2}\right|_{\mathcal L(C([0,1]))}\le C t^{-1/2}. 
$$
This can be seen from the formula
$$
Ae^{-tA^2/2}u= -  \sum_{i\in\N} \lambda_i e^{-\lambda_i^2 t/2} \langle u,e_i \rangle e_i,
$$
where $(e_i)_{i\in\N}$ and $(\lambda_i)_{i\in\N}$ are the eigenvectors and eigenvalues of $-A$. 
Since $|e_i|_{C([0,1]} $ are equibounded, we deduce
$$
\begin{array}{ll}
\ds \left|Ae^{-tA^2/2}u\right|_{C([0,1]} & \ds 
\le C \left(\sum_{i\in\N} \lambda_i e^{-\lambda_i^2 t/2}\right) |u|_{L^1(0,1)}\\
\\
&\le C \, t^{-1/2} |u|_{C([0,1])} 
\end{array}
$$
We deduce 
$$
\begin{array}{l}
\left| Y^{z_1}(t)-Y^{z_2}(t)\right|_{C([0,1])}\\
\\
\le C\, L_\delta \int_0^t (t-s)^{-1/2} 
\left(\left| Y^{z_1}(s)-Y^{z_2}(s)\right|_{C([0,1])} + \left|z_1(s)-z_2(s)\right|_{C([0,1])}
\right) \dds
\end{array}
$$
where $L_\delta$ is the Lipschitz constant of $f_\delta$. Gronwall Lemma implies the result for $T$ 
sufficiently small. Iterating the argument we obtain the continuity of $\Phi_\delta$.
\EProof
\begin{Coro}
\label{c3.1}
For every $c\in (-1,1)$, $\nu_c$ is the unique invariant measure of the transition semigroup
$(P_t^c)_{t\ge 0}$. Moreover it is ergodic.
\end{Coro}
Using classical arguments, it is easily seen that, for $\lambda=0$, $\nu_c^n$ satisfies a log-Sobolev inequality and therefore a Poincar\'e inequality. The constant in these inequality do not depend
on $n$ so that we have the same result for $\nu_c$. For $\lambda\ne 0$, we can argue as in \cite{MR1918538}
 and prove that this is still true.

We now want to prove a stronger result : exponential mixing. We use coupling arguments developped by Odasso in \cite{MR2276442}.
\begin{Th}\label{Th:Mixing}
For every $c \in (-1,1)$, there exist a small $\beta >0$ and a constant $C>0$ such that for all $\varphi\in \mathcal B_b(K\cap H_c)$, $t >0$ and $x\in H_{c}$
\BE\label{Eq:Mixing}
\left|\E[\varphi(X(t,x))] - \nu_{c}(\varphi)\right| \leq C \|\varphi\|_{\infty} e^{-\beta t}.
\EE
\end{Th}
\Proof
By \eqref{Eq:2.8}, we know that for any $\varphi\in \mathcal B_b(K\cap H_c)$, $T> 0$, $\varepsilon >0$,
$$
| P_T^c\varphi(x)- P_T^c\varphi(y)|\le \frac{4e^{\lambda^2T/4}}{\lambda\sqrt{T}} \varepsilon \|\varphi\|_\infty
$$
if $x,\, y\in H_c$, $|x|_{-1}\le \varepsilon$ and  $|y|_{-1}\le \varepsilon$. By definition of the total variation norm, we deduce
\BE
\label{e3.3}
\|\left(P_T^c\right)^*\delta_x - \left(P_T^c\right)^*\delta_y\|_{var}
= \sup_{\|\varphi\|_\infty \le 1} | P_T^c\varphi(x)- P_T^c\varphi(y)|
\le \frac{4e^{\lambda^2T/4}}{\lambda\sqrt{T}} \varepsilon
\EE
for $T>0$, $x,\, y\in H_c$, $|x|_{-1}\le \varepsilon$ and  $|y|_{-1}\le \varepsilon$. We have denoted by 
$\delta_x$ the Dirac mass at $x\in H_c$ so that $\left(P_T^c\right)^*\delta_x$ is the law of $X(T,x)$.

Recall that a coupling of $(\left(P_T^c\right)^*\delta_x,\left(P_T^c\right)^*\delta_y)$ is a couple of random
variable $(X_1,X_2)$ such that the law of $X_1$ is $\left(P_T^c\right)^*\delta_x$ and the law of 
$X_2$ is $\left(P_T^c\right)^*\delta_y$. By standard results on couplings (see for instance \cite{MR1845328} 
section 4, or \cite{MR1937652}), we know there exists a maximal coupling of 
$(\left(P_T^c\right)^*\delta_x,\left(P_T^c\right)^*\delta_y)$. Let us denote by $(Y_1(x,y),Y_2(x,y))$ this
maximal coupling, it satisfies
\BE
\label{e3.3bis}
\P(Y_1(x,y)\ne Y_2(x,y)) = \|\left(P_T^c\right)^*\delta_x - \left(P_T^c\right)^*\delta_y\|_{var}.
\EE
Moreover $(Y_1(x,y),Y_2(x,y))$ depends measurably on $(x,y)$. 

By the Strong Feller property, we know that $x\mapsto \P(|X(T,x)|_{-1}\le \varepsilon)$ is
continuous on $H_c$. Therefore, thanks to Proposition \ref{p3.1}, for any $x\in K\cap H_c$, there exists 
a $\eta_x>0$ and a $\kappa_x>0$ such that
$$
\P(|X(T,y)|_{-1}\le \varepsilon) >\kappa_x
$$
for all $y\in K\cap H_c$ such that $|x-y|_{-1}\le \eta_x$. 
By compactness of $K\cap H_c$ in $H_c$, we 
deduce that that there exits $\kappa_0>0$ such that 
\BE
\label{e3.4}
\P(|X(T,y)|_{-1}\le \varepsilon) >\kappa_0
\EE
for all $y\in K\cap H_c$.

Let $\tilde W$ a cylindrical Wiener process independent on $W$ and denote 
by $\tilde X$ the associated solution of the stochastic Cahn-Hilliard equation which has the same law as $X$. 
For arbitrary $x,\, y \in K\cap H_c$, we define the coupling 
$(Z_1(x,y),Z_2(x,y)$ of 
$(\left(P_T^c\right)^*\delta_x,\left(P_T^c\right)^*\delta_y)$  as follows
$$
(Z_1(x,y),Z_2(x,y))=
\left\{
\begin{array}{l}
(X(T,x), X(T,y)) \mbox{ if } x=y,\\
\\
(Y_1(x,y),Y_2(x,y))  \mbox{ if } |x|_{-1}\le \varepsilon,\; |y|_{-1}\le \varepsilon \mbox{ and } x\ne y,\\
\\
(X(T,x), \tilde X(T,y)) \mbox{ otherwise.}
\end{array}
\right.
$$

We now construct recursively $(X_1(kT,x,y),X_2(kT,x,y))$ a coupling of $(\left(P_{kT}^c\right)^*\delta_x,\left(P_{kT}^c\right)^*\delta_y)$,
the laws of  $X(kT,x)$ and  $X(kT,y)$. For $k=0$, we set $(X_1(kT,x,y),X_2(kT,x,y))=(x,y)$. For 
$k\ge 0$, we define $\left(X_1\left((k+1)T,x,y\right),X_2\left((k+1)T,x,y\right)\right)$ by
\BEAS
X_1\left((k+1)T,x,y\right) &=& Z_1\left(X_1\left(kT,x,y\right),X_2\left(kT,x,y\right)\right),\\
X_2\left((k+1)T,x,y\right) &=& Z_2\left(X_1\left(kT,x,y\right),X_2\left(kT,x,y\right)\right).
\EEAS
Let us define 
$$
\tau= \inf\{ kT : \; |X_1(kT,x,y)|_{-1}\le \varepsilon,\; |X_2(kT,x,y)|_{-1}\le \varepsilon\}
$$

If $|x|_{-1}\le \varepsilon$ and $|y|_{-1}\le \varepsilon$, then $\tau=0$ 
and $\E(e^{\alpha \tau})=1$.

If $\tau\ne 0$ {\it i.e.} if $|x|_{-1}\ge \varepsilon$ or $|y|_{-1}\ge \varepsilon$, then by construction of the coupling and \eqref{e3.4}
$$
\P(\tau >T)< 1-Ê\kappa_0^2.
$$
More generally
$$
\P(\tau >kT\big| \tau\ge kT) <  1-Ê\kappa_0^2.
$$
We deduce
$$
\P(\tau >kT) <  (1-Ê\kappa_0^2)^k
$$
and 
$$
\E(e^{\alpha \tau})=\sum_{k\in \N} e^{\alpha kT}\P(\tau=kT)\le \sum_{k\in \N} e^{\alpha kT}
(1-Ê\kappa_0^2)^{k-1} =M<\infty
$$
for $\alpha$ small enough.
Similarly, if we define 
$$
\tau_n= \inf\{ kT>\tau_{n-1} : \; |X_1(kT,x,y)|_{-1}\le \varepsilon,\; |X_2(kT,x,y)|_{-1}\le \varepsilon\},
$$
for all $n \geq 2$ and with $\tau_{1} := \tau$.
We have 
$$
\E(e^{\alpha (\tau_n-\tau_{n-1})})\le M
$$
so that
$$
\E(e^{\alpha \tau_n})\le M^n.
$$
Define 
$$
k_0=\inf\{n\geq1 :\; X_1(\tau_n+T,x,y)=X_2(\tau_n+T,x,y)\}.
$$
By \eqref{e3.3}, \eqref{e3.3bis}, for all $n \geq 1$
$$
\P(k_0=n)\leq \P(k_0>n-1)\leq \left(\frac{4e^{\lambda^2T/4}}{\lambda\sqrt{T}} \varepsilon\right)^{n-1}.
$$
We choose $\varepsilon$ small enough such that
\BD
\left(\frac{4e^{\lambda^2T/4}}{\lambda\sqrt{T}} \varepsilon\right) <1.
\ED
Then we write
\BEAS
\E(e^{\beta \tau_{k_0}})=\sum_{n\geq 1} \E\left(e^{\beta \tau_n} \1_{k_0=n}\right)&\leq&
\sum_{n\geq 1} (\E(e^{2 \beta \tau_n}))^{1/2} (\P(k_0=n))^{1/2}\\
&\leq& \sum_{n\geq 1} M^{n\beta/\alpha}
\left(\frac{4e^{\lambda^2T/4}}{\lambda\sqrt{T}} \varepsilon\right)^{(n-1)/2}=\mathcal{M}<\infty
\EEAS
for $\beta$ small enough.

By Markov's inequality, we conclude that for all $k\geq 1$
$$
\begin{array}{l}
|\E(\varphi(X(kT,x)))-\E(\varphi(X(kT,y)))| \\
\\
=|\E(\varphi(X_1(kT,x,y)))-\E(\varphi(X_2(kT,x,y)))|\\
\\
\le 2\|\varphi\|_\infty \P(X_1(kT,x,y)\ne X_2(kT,x,y))\\
\\
\le 2 \|\varphi\|_\infty \P(kT>\tau_{k_0}+T)  \\
\\
\le 2 \|\varphi\|_\infty\mathcal{M} e^{-\beta(k-1)T}.
\end{array}
$$
We define $k:= \left\lfloor \frac{t}{T}\right\rfloor$ such that we have $P_{t}^c = P_{kT}^cP_{t-kT}^c$. Thus we can write
\BEAS
\left|\E[\varphi(X(t,x))] - \nu_{c}(\varphi)\right| &=& \left|P_{t}^c\varphi(x) - \int_{H_{c}} \varphi(y) \nu_{c}(\ddy)\right| \\
&=& \left|\int_{H_{c}}P_{t}^c\varphi(x)\nu_{c}(\ddy) - \int_{H_{c}} P_{t}^c\varphi(y) \nu_{c}(\ddy)\right| \\
&\leq &\left|\int_{H_{c}} \left(P_{kT}^cP_{t-kT}^c\varphi(x) - P_{kT}^cP_{t-kT}^c\varphi(y)\right) \nu_{c}(\ddy)\right|\\
&\leq &\int_{H_{c}} 2  \|P_{t-kT}^c\varphi\|_\infty\mathcal{M} e^{-\beta(k-1)T} \nu_{c}(\ddy)\\
&\leq & 2  \|\varphi\|_\infty\mathcal{M} e^{-\beta(k-1)T}\\
&\leq & C \|\varphi\|_\infty e^{-\beta t}.
\EEAS
\EProof

\bibliographystyle{abbrv}

\bibliography{Bibliludo}
\end{document}